\documentclass[a4paper,11pt]{amsart}
\usepackage{graphicx}
\usepackage{mathptmx}
\usepackage{amsmath}
\usepackage{amssymb}
\usepackage{enumitem}
\usepackage{xcolor}

\newmuskip\pFqmuskip

\newcommand*\pFq[6][8]{%
  \begingroup 
  \pFqmuskip=#1mu\relax
  \mathcode`=\string"8000
  \begingroup\lccode`\~=`\,
  \lowercase{\endgroup\let~}\pFqcomma
  F^{#2}_{#3}{\left(\genfrac..{0pt}{}{#4}{#5}\bigg|#6\right)}%
  \endgroup
}
\newcommand{\pFqcomma}{\mskip\pFqmuskip}

\newtheorem{theorem}{Theorem}

\newtheorem{proposition}[theorem]{Proposition}

\begin{document}

\title[]{Some identities on degenerate hyperbolic functions arising from $p$-adic integrals on $\mathbb{Z}_p$}

\author{Taekyun  Kim}
\address{Department of Mathematics, Kwangwoon University, Seoul 139-701, Republic of Korea}
\email{tkkim@kw.ac.kr}

\author{Hye Kyung  Kim}
\address{Department of Mathematics Education, Daegu Catholic University, Gyeongsan 38430, Republic of Korea}
\email{hkkim@cu.ac.kr}

\author{Dae San Kim}
\address{Department of Mathematics, Sogang University, Seoul 121-742, Republic of Korea}
\email{dskim@sogang.ac.kr}

\thanks{* is corresponding author}

\begin{abstract}
The aim of this paper is to introduce several degenerate hyperbolic functions as degenerate versions of the hyperbolic functions, to evaluate Volkenborn and the fermionic $p$-adic integrals of the degenerate hyperbolic cosine and the degenerate hyperbolic sine functions and to derive from them some identities involving the degenerate Bernoulli numbers, the degenerate Euler numbers and the Cauchy numbers of the first kind.
\end{abstract}

 \maketitle

\markboth{\centerline{\scriptsize  Some identities of degenerate hyperbolic functions arising from $p$-adic integrals on $\mathbb{Z}_p$}}
{\centerline{\scriptsize  T. Kim, H. K. Kim and D. S. Kim}}


\section{Introduction}

Volkenborn and the fermionic $p$-adic integrals of the powers of $x$ give respectively the Bernoulli numbers and the Euler numbers, while Volkenborn and the fermionic $p$-adic integrals of the generalized falling factorials yield respectively the degenerate Bernoulli numbers and the degenerate Euler numbers (see \eqref{eq10}). Thus the latter may be viewed as degenerate versions of the former. \par
The aim of this paper is to derive some identities on degenerate hyperbolic functions arising from Volkenborn and the fermionic $p$-adic integrals on $\mathbb{Z}_{p}$. In more detail, we introduce degenerate hyperbolic functions as natural degenerate versions of the usual hyperbolic functions and derive some identities relating to them. Then we evaluate Volkenborn and the fermionic $p$-adic integrals of the degenerate hyperbolic cosine and the degenerate hyperbolic sine functions. From those results, we derive some identities involving the degenerate Bernoulli numbers, the degenerate Euler numbers and the Cauchy numbers of the first kind.  \par
The outline of this paper is as follows. In Section 1, we recall Volkenborn integral of uniformly differentiable functions on $\mathbb{Z}_{p}$ and the fermionic $p$-adic  integral of continuous functions on $\mathbb{Z}_{p}$, together with their integral equations and examples in the case of exponential functions. We remind the reader of Cauchy numbers of the first kind. We recall the degenerate exponentials as a degenerate version of the usual exponentials. Then we show that the degenerate Bernoulli and the degenerate Euler numbers arise naturally respectively from the Volkenborn and the fermionic $p$-adic integrals on $\mathbb{Z}_{p}$ of the degenerate exponentials. Section 2 is the main result of this paper. We introduce degenerate versions of the hyperbolic functions, namely the degenerate hyperbolic cosine $\cosh_{\lambda}(x:a)$, the degenerate hyperbolic sine $\sinh_{\lambda}(x:a)$, degenerate hyperbolic tangent $\tanh_{\lambda}(x:a)$ and the degenerate hyperbolic cotangent $\coth_{\lambda}(x:a)$. Then we derive several identities connecting  those degenerate hyperbolic functions. We evaluate Volkenborn integrals of $\cosh_{\lambda}(x:a)$ and $\sinh_{\lambda}(x:a)$. From those results, we obtain some identities involving the degenerate Bernoulli numbers and the Cauchy numbers of the first kind. We compute the fermionic $p$-adic integrals of the same degenerate hyperbolic functions, from which we derive an identity involving the degenerate Euler numbers. Finally, we conclude this paper in Section 3. In the rest of this section, we recall the necessary facts that are needed throughout this paper.\par

Let $p$ be a fixed odd prime mber. Throughout this paper, $\mathbb{Z}_p, \mathbb{Q}_p$ and $\mathbb{C}_p$ will denote the ring of $p$-adic integers, the field of $p$-adic rational numbers and the completion of an algebraic closure of $\mathbb{Q}_p$, respectively. The $p$-adic norm $|\cdot|_p$ is normalized as $|p|_p=\frac{1}{p}$. The Volkenborn integral on $\mathbb{Z}_p$ is defined by


\begin{equation}\label{eq01}
\begin{split}
\int_{\mathbb{Z}_p}f(x) \ d\mu(x)&=\lim_{N \rightarrow\infty}\sum_{x=0}^{P^N-1}f(x)\mu(x+p^N\mathbb{Z}_p)\\
&=\lim_{N\rightarrow\infty}\frac{1}{p^N}\sum_{x=0}^{p^N-1}f(x), \quad (\text{see \cite{1,2,5,6,7,8,13}}).
\end{split}
\end{equation}

From \eqref{eq01}, we note that
\begin{equation}\label{eq02}
\begin{split}
\int_{\mathbb{Z}_p}f(x+1) \ d\mu(x)-\int_{\mathbb{Z}_p}f(x) \ d\mu(x)=f'(0), \quad (\text{see \cite{5,13}}),
\end{split}
\end{equation}
where $f$ is a $\mathbb{C}_p$-valued uniform differentiable function on $\mathbb{Z}_p$.

In \cite{7}, the $fermionic$ $p$-adic integral on $\mathbb{Z}_p$ is defined by
\begin{equation}\label{eq03}
\begin{split}
\int_{\mathbb{Z}_p}f(x) \ d\mu_{-1}(x)&=\lim_{N\rightarrow\infty}\sum_{x=0}^{p^N-1}f(x) \ \mu_{-1}(x+p^N\mathbb{Z}_p)\\
&=\lim_{N\rightarrow\infty}\sum_{x=0}^{p^N-1}f(x)(-1)^x,
\end{split}
\end{equation}
where $f$ is a $\mathbb{C}_p$-valued continuous function on $\mathbb{Z}_p$.

Thus, by \eqref{eq03}, we get
\begin{equation}\label{eq04}
\begin{split}
\int_{\mathbb{Z}_p}f(x+1) \ d\mu_{-1}(x)+\int_{\mathbb{Z}_p}f(x) \ d\mu_{-1}(x)=2f(0) \quad (\text{see \cite{1,7,8}}).
\end{split}
\end{equation}

Adopting the notation from \cite{13}, we let $E$ denote the the additive group given by
\begin{equation}
\begin{split}\label{eq04-1}
E=\left\{x \in \mathbb{C}_p\, \big\vert \,  |x|_{p}<p^{-\frac{1}{p-1}} \right\}.
\end{split}
\end{equation}
It is known that there is an isomorphism from the multiplicative group $1+E$ to the additive group $E$ given by $1+x \mapsto \log(1+x)$, where $\log =\log_{p} $ is the $p$-adic logarithm. Moreover, the inverse map is given by $x \mapsto e^{x}=\exp x$, where $\exp = \exp_{p} $ is the $p$-adic exponential (see \cite{13}).

The Cauchy numbers of the first kind (also called Bernoulli numbers of the second) are defined by 
\begin{equation*}
\begin{split}
\frac{t}{\log(1+t)}=\sum_{n=0}^\infty C_n\frac{t^n}{n!},\quad (|t|_{p}<1,\, t \ne 0).
\end{split}
\end{equation*}

From \eqref{eq02} and \eqref{eq04}, we note that
\begin{equation}\label{eq05}
\begin{split}
\int_{\mathbb{Z}_p}e^{xt} \ d\mu(x)=\frac{t}{e^t-1}=\sum_{n=0}^\infty B_n\frac{t^n}{n!}, \quad (t \in E, t \ne 0), \quad (\text{see [1-18]}),
\end{split}
\end{equation}

 and
\begin{equation}\label{eq06}
\begin{split}
\int_{\mathbb{Z}_p}e^{xt} \ d\mu_{-1}(x)=\frac{2}{e^t+1}=\sum_{n=0}^\infty E_n\frac{t^n}{n!}, \quad (t \in E ),\quad (\text{see \cite{1,5,7}}),
\end{split}
\end{equation}
where $B_n$ and $E_n$ are the $n$-th Bernoulli number and the $n$-th Euler number, respectively.

{\it{Unless otherwise stated, throughout this paper, we let $\lambda$ be any nonzero element in $\mathbb{C}_{p}$.}} \\
The degenerate exponentials are defined by
\begin{equation}\label{eq07}
\begin{split}
e^x_\lambda(t)=\sum_{n=0}^\infty (x)_{n,\lambda}\frac{t^n}{n!}=(1+\lambda t)^{\frac{x}{\lambda}},\quad \big(t \in \frac{1}{\lambda}E,\, |x|_{p} \le |\lambda|_{p}\big), \quad (\text{see \cite{9,10}}),
\end{split}
\end{equation}
where the generlaized falling factorials (also called the $\lambda$-falling factorials) are given by $(x)_{0,\lambda}=1, (x)_{n,\lambda}=x(x-\lambda) \cdots (x-(n-1)\lambda), \quad (n\geq1)$.

From \eqref{eq02}, we observe that
\begin{equation}\label{eq08}
\begin{split}
\int_{\mathbb{Z}_p}e^{x}_\lambda(t) \ d\mu(x)=\frac{\frac{1}{\lambda}\log (1+\lambda t)}{e_\lambda(t)-1}=\sum_{n=0}^\infty \beta_{n,\lambda}\frac{t^n}{n!},\quad \big(t \in \frac{1}{\lambda}E, \, t \ne 0, \, |\lambda|_{p} \ge 1 \big),
\end{split}
\end{equation}
where $\beta_{n,\lambda}$ are called the (fully) degenerate Bernoulli numbers $(\text{see \cite{3,5}})$. \\
Note that $\lim_{\lambda \rightarrow 0}\beta_{n,\lambda}=B_n, \ (n\geq0), \quad (\text{see \cite{11}})$.

From \eqref{eq04}, we note that
\begin{equation}\label{eq09}
\begin{split}
\int_{\mathbb{Z}_p}e^{x}_\lambda(t) \ d\mu_{-1}(x)=\frac{2}{e_{\lambda}(t)+1}=\sum_{n=0}^\infty \mathcal{E}_{n,\lambda}\frac{t^n}{n!},\quad \big(t \in \frac{1}{\lambda}E, \, |\lambda|_{p} \ge 1 \big),
\end{split}
\end{equation}
where $\mathcal{E}_{n,\lambda}$ are called the degenerate Euler numbers $(\text{see \cite{3,5}})$. \\
Note that $\lim_{\lambda\rightarrow0}\mathcal{E}_{n,\lambda}=E_n, \ (n\geq0)$. \\
Thus, by \eqref{eq08} and \eqref{eq09}, we get
\begin{equation}\label{eq10}
\begin{split}
\int_{\mathbb{Z}_p}(x)_{n,\lambda} \ d\mu(x)=\beta_{n,\lambda}, \ \ \int_{\mathbb{Z}_p}(x)_{n,\lambda} \ d\mu_{-1}(x)=\mathcal{E}_{n,\lambda}, \ (n\geq0).
\end{split}
\end{equation}

The hyperbolic functions of real or complex variables are given by
\begin{equation}\label{eq11}
\begin{split}
&\cosh{x}=\frac{e^x+e^{-x}}{2}, \quad \sinh{x}=\frac{e^x-e^{-x}}{2}, \\
&\tanh{x}=\frac{\sinh{x}}{\cosh{x}}, \quad {\rm{and}} \quad \coth{x}=\frac{\cosh{x}}{\sinh{x}}.
\end{split}
\end{equation}

\bigskip

\section{Some identities of degenerate hyperbolic functions arising from $p$-adic integrals on $\mathbb{Z}_p$}
Let $\lambda$ be any nonzero element in $\mathbb{C}_{p}$, and let $ a \in \frac{1}{\lambda}E $,\,(see \eqref{eq04-1}).
Let us consider the degenerate hyperbolic functions given by
\begin{equation}\label{eq12}
\begin{split}
&\cosh_{\lambda}(x:a)=\frac{1}{2}(e^x_\lambda(a)+e^{-x}_\lambda(a)), \quad (|x|_{p} \le |\lambda|_{p}),\\
&\sinh_{\lambda}(x:a)=\frac{1}{2}(e^x_\lambda(a)-e^{-x}_\lambda(a)), \quad (|x|_{p} \le |\lambda|_{p}),\\
&\tanh_{\lambda}(x:a)=\frac{\sinh_{\lambda}(x:a)}{\cosh_{\lambda}(x:a)},\quad (|x|_{p} \le |\lambda|_{p}),\\ 
&\coth_{\lambda}(x:a)=\frac{\cosh_{\lambda}(x:a)}{\sinh_{\lambda}(x:a)},\quad (|x|_{p} \le |\lambda|_{p}, \,x \ne 0, \,a \ne 0).
\end{split}
\end{equation}
Note that $\lim_{\lambda\rightarrow0}\cosh_{\lambda}(x:a)=\cosh(xa), \ \  \lim_{\lambda\rightarrow0}\sinh_{\lambda}(x:a)=\sinh(xa), \\ \lim_{\lambda\rightarrow0}\tanh_{\lambda}(x:a)=\tanh(xa), \ \ {\rm{and}} \ \ \lim_{\lambda\rightarrow0}\coth_{\lambda}(x:a)=\coth(xa).$\\

From \eqref{eq12}, we note that
\begin{equation}\label{eq13}
\begin{split}
\cosh_{\lambda}(2x:a)=&\frac{e^{2x}_\lambda(a)+e^{-2x}_\lambda(a)}{2}=1+\frac{e^{2x}_\lambda(a)+e^{-2x}_\lambda(a)-2}{2}\\
&=2\bigg(\frac{e^x_\lambda(a)-e^{-x}_\lambda(a)}{2}\bigg)^2+1=1+\sinh^2_\lambda(x:a), \quad (|x|_{p} \le |\lambda|_{p}).
\end{split}
\end{equation}

On the other hand, by \eqref{eq12}, we get
\begin{equation}\label{eq14}
\begin{split}
2\cosh^2_\lambda(x:a)-1&=2\times\frac{e^{2x}_\lambda(a)+e^{-2x}_\lambda(a)+2}{4}-1\\
&=\frac{e^{2x}_\lambda(a)+e^{-2x}_\lambda(a)}{2}=\cosh_\lambda(2x:a), \quad (|x|_{p} \le |\lambda|_{p}).
\end{split}
\end{equation}

Therefore, by \eqref{eq13} and \eqref{eq14}, we obtain the following proposition.
\begin{proposition}
Let $a \in \frac{1}{\lambda}E$, and let $|x|_{p} \le |\lambda|_{p} $. Then the following identity holds true.
\begin{equation*}
\begin{split}
\cosh_\lambda(2x:a)&=2\cosh^2_\lambda(x:a)-1=1+2\sinh^2_\lambda(x:a).
\end{split}
\end{equation*}
\end{proposition}

\medskip

By \eqref{eq12}, we get
\begin{equation}\label{eq15}
\begin{split}
\sinh_\lambda(2x:a)&=\frac{e^{2x}_\lambda(a)-e^{-2x}_\lambda(a)}{2}\\
&=2\frac{e^{x}_\lambda(a)-e^{-x}_\lambda(a)}{2}\frac{e^{x}_\lambda(a)+e^{-x}_\lambda(a)}{2}\\
&=2\sinh_\lambda(x:a)\cosh_\lambda(x:a).
\end{split}
\end{equation}

\medskip

Thus we have the following proposition.
\begin{proposition}Let $a \in \frac{1}{\lambda}E$, and let $|x|_{p} \le |\lambda|_{p} $. Then the following identity holds true.
\begin{equation*}
\begin{split}
\sinh_\lambda(2x:a)=2\sinh_\lambda(x:a)\cosh_\lambda(x:a).
\end{split}
\end{equation*}
\end{proposition}

Now, we observe that
\begin{equation}\label{eq16}
\begin{split}
&\sinh_\lambda(x+y:a)=\frac{1}{2}(e_\lambda^{x+y}(a)-e_\lambda^{-x-y}(a))\\
&=\frac{e_\lambda^{x+y}(a)-e_\lambda^{-x-y}(a)+e_\lambda^{x-y}(a)-e_\lambda^{-x+y}(a)}{4}+\frac{e_\lambda^{x+y}(a)-e_\lambda^{-x-y}(a)-e_\lambda^{x-y}(a)+e_\lambda^{-x+y}(a)}{4}\\
&=\frac{e_\lambda^{x}(a)-e_\lambda^{-x}(a)}{2}\times\frac{e_\lambda^{y}(a)+e_\lambda^{-y}(a)}{2}+\frac{e_\lambda^{x}(a)+e_\lambda^{-x}(a)}{2}\times\frac{e_\lambda^{y}(a)-e_\lambda^{-y}(a)}{2}\\
&=\sinh_\lambda(x:a)\cosh_\lambda(y:a)+\cosh_\lambda(x:a)\sinh_\lambda(y:a),\quad (a \in \frac{1}{\lambda}E,\,|x|_{p} \le |\lambda|_{p},\,|y|_{p} \le |\lambda|_{p}).
\end{split}
\end{equation}

On the other hand, by \eqref{eq12}, we get
\begin{equation}\label{eq17}
\begin{split}
&\cosh_\lambda(x:a)\cosh_\lambda(y:a)+\sinh_\lambda(x:a)\sinh_\lambda(y:a)\\
&=\frac{e_\lambda^{x}(a)+e_\lambda^{-x}(a)}{2}\times\frac{e_\lambda^{y}(a)+e_\lambda^{-y}(a)}{2}+\frac{e_\lambda^{x}(a)-e_\lambda^{-x}(a)}{2}\times\frac{e_\lambda^{y}(a)-e_\lambda^{-y}(a)}{2}\\
&=\frac{e_\lambda^{x+y}(a)+e_\lambda^{-x-y}(a)+e_\lambda^{-x+y}(a)+e_\lambda^{x-y}(a)}{4}+\frac{e_\lambda^{x+y}(a)+e_\lambda^{-x-y}(a)-e_\lambda^{x-y}(a)-e_\lambda^{-x+y}(a)}{4}\\
&=\frac{2(e_\lambda^{x+y}(a)+e_\lambda^{-(x+y)}(a))}{4}=\frac{e_\lambda^{x+y}(a)+e_\lambda^{-(x+y)}(a)}{2}=\cosh_\lambda(x+y:a),
\end{split}
\end{equation}
where $(a \in \frac{1}{\lambda}E,\,|x|_{p} \le |\lambda|_{p},\,|y|_{p} \le |\lambda|_{p})$.

Therefore, by \eqref{eq16} and \eqref{eq17}, we obtain the following proposition.
\begin{proposition} Let $a \in \frac{1}{\lambda}E,\,|x|_{p} \le |\lambda|_{p}$,\, and let $|y|_{p} \le |\lambda|_{p}$. Then the following identities are valid.
\begin{equation*}
\begin{split}
&\cosh_\lambda(x+y:a)=\cosh_\lambda(x:a)\cosh_\lambda(y:a)+\sinh_\lambda(x:a)\sinh_\lambda(y:a), \\
&{\rm{and}} \quad \sinh_\lambda(x+y:a)=\sinh_\lambda(x:a)\cosh_\lambda(y:a)+\cosh_\lambda(x:a)\sinh_\lambda(y:a).
\end{split}
\end{equation*}
\end{proposition}

\medskip

Noting that $\frac{d}{dx}\cosh_\lambda(x:a)=\frac{1}{\lambda}\log(1+\lambda a)\sinh_\lambda(x:a)$ and using \eqref{eq02}, we have

\begin{equation}\label{eq18}
\begin{split}
0&=\int_{\mathbb{Z}_p}\cosh_\lambda(x+1:a) \ d\mu(x)-\int_{\mathbb{Z}_p}\cosh_\lambda(x:a) \ d\mu(x)\\
&=(\cosh_\lambda(1:a)-1)\int_{\mathbb{Z}_p}\cosh_\lambda(x:a) \ d\mu(x)+\sinh_\lambda(1:a)\int_{\mathbb{Z}_p}\sinh_\lambda(x:a) \ d\mu(x),
\end{split}
\end{equation}
where $\big(a \in \frac{1}{\lambda}E,\,|\lambda|_{p} \ge 1 \big)$.

Observing that $\frac{d}{dx}\sinh_\lambda(x:a)=\frac{1}{\lambda}\log{(1+\lambda a)}\cosh_\lambda(x:a)$ and making use of \eqref{eq02}, we get

From \eqref{eq02} we note that
\begin{equation}\label{eq20}
\begin{split}
\frac{1}{\lambda}\log(1+\lambda a)
&=\int_{\mathbb{Z}_p}\sinh_{\lambda}(x+1:a) d\mu(x)-\int_{\mathbb{Z}_p}\sinh_{\lambda}(x:a) d\mu(x)\\
&=(\cosh_\lambda(1:a)-1)\int_{\mathbb{Z}_p}\sinh_\lambda(x:a) \ d\mu(x) \\
&\quad\quad +\sinh_\lambda(1:a)\int_{\mathbb{Z}_p}\cosh_\lambda(x:a) \ d\mu(x),
\end{split}
\end{equation}
where $\big(a \in \frac{1}{\lambda}E,\,|\lambda|_{p} \ge 1 \big)$.

By solving the system of linear equations in \eqref{eq18} and \eqref{eq20}, we get
\begin{equation}\label{eq21}
\begin{split}
&\int_{\mathbb{Z}_p}\cosh_\lambda(x:a) \ d\mu(x)
=\frac{1}{\lambda}\log(1+\lambda a) \frac{\sinh_\lambda(1:a)}{2(\cosh_\lambda(1:a)-1)}, \\
&\int_{\mathbb{Z}_p}\sinh_\lambda(x:a) \ d\mu(x)=-\frac{1}{2}\frac{1}{\lambda}\log(1+\lambda a),\quad \big(a \in \frac{1}{\lambda}E,\,|\lambda|_{p} \ge 1 \big).
\end{split}
\end{equation}

From \eqref{eq21} and Propositions 1 and 2, we have
\begin{equation}\label{eq23}
\begin{split}
\int_{\mathbb{Z}_p}\cosh_\lambda(x:a)d\mu(x)&=\frac{1}{\lambda}\log(1+\lambda a)
\frac{2\sinh_\lambda(\frac{1}{2}:a)\cosh_\lambda(\frac{1}{2}:a)}{4\sinh_\lambda^2(\frac{1}{2}:a)}\\
&=\frac{a}{2}\coth_\lambda\bigg(\frac{1}{2}:a\bigg)\frac{\log(1+\lambda a)}{\lambda a},
\quad \big(a \in \frac{1}{\lambda}E,\,|\lambda|_{p} \ge 1 \big).
\end{split}
\end{equation}

Therefore, from \eqref{eq21} and \eqref{eq23} we obtain the following theorem.
\begin{theorem} Let $ a \in \frac{1}{\lambda}E$,\,and let $|\lambda|_{p} \ge 1 $. Then the following relations hold true.
\begin{equation*}
\begin{split}
&\frac{\lambda a}{\log(1+\lambda a)}\int_{\mathbb{Z}_p}\cosh_\lambda(x:a) \ d\mu(x)=\frac{a}{2}\coth_\lambda\bigg(\frac{1}{2}:a\bigg),\\
&\frac{\lambda a}{\log(1+\lambda a)}\int_{\mathbb{Z}_p}\sinh_\lambda(x:a) \ d\mu(x)
=-\frac{a}{2}.
\end{split}
\end{equation*}
\end{theorem}

\medskip

By Taylor expansion, we get
\begin{equation}\label{eq24}
\begin{split}
\cosh_\lambda(x:a)&=\frac{1}{2}(e^x_\lambda(a)+e^{-x}_\lambda(a))=\frac{1}{2}(e^x_\lambda(a)+e^{x}_{-\lambda}(-a))\\
&=\frac{1}{2}\sum_{m=0}^\infty((x)_{m,\lambda}+(-1)^m(x)_{m,-\lambda})\frac{a^m}{m!},\quad (a \in \frac{1}{\lambda}E,\,|x|_{p} \le |\lambda|_{p}).
\end{split}
\end{equation}

From Theorem 4, \eqref{eq24} and \eqref{eq10}, we have
\begin{equation}\label{eq25}
\begin{split}
\frac{a}{2}\coth_\lambda\bigg(\frac{1}{2}:a\bigg)&=\frac{\lambda}{\log(1+\lambda a)}\int_{\mathbb{Z}_p}\cosh_\lambda(x:a) \ d\mu(x)\\
&=\bigg(\sum_{l=0}^\infty{\lambda}^l C_l\frac{a^l}{l!}\bigg)\frac{1}{2}\sum_{m=0}^\infty\int_{\mathbb{Z}_p}((x)_{m,\lambda}+(-1)^m(x)_{m,-\lambda})d\mu(x)\frac{a^m}{m!}\\
&=\sum_{l=0}^\infty{\lambda}^l C_l\frac{a^l}{l!}\sum_{m=0}^\infty\bigg(\frac{\beta_{m,\lambda}+(-1)^m\beta_{m,-\lambda}}{2}\bigg)\frac{a^m}{m!}\\
&=\sum_{n=0}^\infty\bigg(\sum_{m=0}^n\binom{n}{m}\bigg(\frac{\beta_{m,\lambda}+(-1)^m\beta_{m,-\lambda}}{2}\bigg)C_{n-m}\lambda^{n-m}\bigg)\frac{a^n}{n!},
\end{split}
\end{equation}
where $a \in \frac{1}{\lambda}E$,\,and let $|\lambda|_{p} \ge 1$.
By Taylor expansion, we get
\begin{equation}\label{eq27}
\begin{split}
\sinh_\lambda(x:a)=\frac{1}{2}\sum_{m=0}^\infty\big((x)_{m,\lambda}-(-1)^m(x)_{m,-\lambda}\big)\frac{a^m}{m!},\quad (a \in \frac{1}{\lambda}E,\,|x|_{p} \le |\lambda|_{p}).
\end{split}
\end{equation}

Thus, by \eqref{eq21} and \eqref{eq27}, we get
\begin{equation}\label{eq28}
\begin{split}
-\frac{1}{2 \lambda}\log(1+\lambda a)&=\int_{\mathbb{Z}_p}\sinh_\lambda(x:a) \ d\mu(x) \\
&=\sum_{n=0}^\infty\bigg(\frac{\beta_{n,\lambda}-(-1)^n\beta_{n,-\lambda}}{2}\bigg)\frac{a^n}{n!}=\sum_{n=1}^\infty\bigg(\frac{\beta_{n,\lambda}-(-1)^n\beta_{n,-\lambda}}{2}\bigg)\frac{a^n}{n!},
\end{split}
\end{equation}
$\big(a \in \frac{1}{\lambda}E,\,|\lambda|_{p} \ge 1 \big)$.

\vspace{0.1in}
\noindent and

\begin{equation}\label{eq29}
\begin{split}
-\frac{1}{2\lambda}\log(1+\lambda a)&=-\frac{1}{2}\sum_{n=1}^\infty \frac{(-1)^{n-1}\lambda^{n-1}}{n}a^n=-\frac{1}{2}\sum_{n=1}^\infty(n-1)!(-1)^{n-1}\lambda^{n-1}\frac{a^n}{n!},\quad (|\lambda a|_{p} \le 1).
\end{split}
\end{equation}

By \eqref{eq25}, \eqref{eq28} and \eqref{eq29}, we obtain the following theorem.

\begin{theorem}
For $n \in \mathbb{N}$, we have the identity:
\begin{equation*}
\begin{split}
-(n-1)!(-\lambda)^{n-1}=\big(\beta_{n,\lambda}-(-1)^{n}\beta_{n,-\lambda}\big).
\end{split}
\end{equation*}
In addition, we have the following relation:
\begin{equation*}
\begin{split}
\frac{a}{2}\coth_\lambda\bigg(\frac{1}{2}:a\bigg)=\sum_{n=0}^\infty\bigg(\sum_{m=0}^n\binom{n}{m}\bigg(\frac{\beta_{m,\lambda}+(-1)^m\beta_{m,-\lambda}}{2}\bigg)C_{n-m}\lambda^{n-m}\bigg)\frac{a^n}{n!},
\end{split}
\end{equation*}
where $a \in \frac{1}{\lambda} E$,\,and $|\lambda|_{p} \ge 1$.
\end{theorem}

\medskip

From \eqref{eq04}, we note that
\begin{equation}\label{eq30}
\begin{split}
2&=\int_{\mathbb{Z}_p}\cosh_\lambda(x+1:a) \ d\mu_{-1}(x)+\int_{\mathbb{Z}_p}\cosh_\lambda(x:a) \ d\mu_{-1}(x) \\
&=(\cosh_\lambda(1:a)+1)\int_{\mathbb{Z}_p}\cosh_\lambda(x:a) \ d\mu_{-1}(x)+\sinh_\lambda(1:a)\int_{\mathbb{Z}_p}\sinh_\lambda(x:a) \ d\mu_{-1}(x),
\end{split}
\end{equation}

and

\begin{equation}\label{eq31}
\begin{split}
0&=\int_{\mathbb{Z}_p}\sinh_\lambda(x+1:a) \ d\mu_{-1}(x)+\int_{\mathbb{Z}_p}\sinh_\lambda(x:a) \ d\mu_{-1}(x) \\
&=(\cosh_\lambda(1:a)+1)\int_{\mathbb{Z}_p}\sinh_\lambda(x:a) \ d\mu_{-1}(x)+\sinh_\lambda(1:a)\int_{\mathbb{Z}_p}\cosh_\lambda(x:a) \ d\mu_{-1}(x),
\end{split}
\end{equation}
where $a \in \frac{1}{\lambda}E$,\,and $|\lambda|_{p} \ge 1$.

By solving the system of linear equations in \eqref{eq30} and \eqref{eq31}, we obtain

\begin{equation}\label{eq32}
\begin{split}
&\int_{\mathbb{Z}_p}\cosh_\lambda(x:a) \ d\mu_{-1}(x)=1, \\
&\int_{\mathbb{Z}_p}\sinh_\lambda(x:a) \ d\mu_{-1}(x)=-\frac{\sinh_{\lambda}(1:a)}{\cosh_{\lambda}(1:a)+1},
\end{split}
\end{equation}
where $a \in \frac{1}{\lambda} E$,\,and $|\lambda|_{p} \ge 1$.

Thus, by \eqref{eq32} and Propositions 1 and 2, we get
\begin{equation}\label{eq34}
\begin{split}
\int_{\mathbb{Z}_p}\sinh_\lambda(x:a) \ d\mu_{-1}(x)=-\frac{2\sinh_\lambda(\frac{1}{2}:a)\cosh_\lambda(\frac{1}{2}:a)}{2\cosh_\lambda^2(\frac{1}{2}:a)}
=-\tanh_\lambda\bigg(\frac{1}{2}:a\bigg),
\end{split}
\end{equation}
where $a \in \frac{1}{\lambda} E$,\,and $|\lambda|_{p} \ge 1$.

Therefore, by \eqref{eq32} and \eqref{eq34}, we obtain the following theorem.

\begin{theorem}
Let $a \in \frac{1}{\lambda} E$,\,and let $|\lambda|_{p} \ge 1$. Then we have the following relations.
\begin{equation*}
\begin{split}
&\int_{\mathbb{Z}_p}\sinh_\lambda(x:a) \ d\mu_{-1}(x)=-\tanh_\lambda\bigg(\frac{1}{2}:a \bigg), \\
& \quad \int_{\mathbb{Z}_p}\cosh_\lambda(x:a) \ d\mu_{-1}(x)=1.
\end{split}
\end{equation*}
\end{theorem}

Thus, by Theorem 6, \eqref{eq27} and \eqref{eq10}, we get
\begin{equation}\label{eq37}
\begin{split}
-\tanh_\lambda(\frac{1}{2}:a)&=\int_{\mathbb{Z}_p}\sinh_\lambda(x:a) \ d\mu_{-1}(x) \\
&=\sum_{n=0}^\infty\bigg(\frac{\mathcal{E}_{n,\lambda}-(-1)^n\mathcal{E}_{n,-\lambda}}{2}\bigg)\frac{a^n}{n!},\quad \big(a \in \frac{1}{\lambda}E,\,|\lambda|_{p} \ge 1 \big).
\end{split}
\end{equation}

In addition, by Theorem 6, \eqref{eq24} and \eqref{eq10}, we get
\begin{equation}\label{eq39}
\begin{split}
1&=\int_{\mathbb{Z}_p}\cosh_\lambda(x;a) \ d\mu_{-1}(x) =\sum_{n=0}^\infty \bigg(\frac{\mathcal{E}_{n,\lambda}+(-1)^n\mathcal{E}_{n,-\lambda}}{2}\bigg)\frac{a^n}{n!},,\quad \big(a \in \frac{1}{\lambda}E,\,|\lambda|_{p} \ge 1 \big).
\end{split}
\end{equation}

Thus, by\eqref{eq39}, we get
\begin{equation}\label{eq40}
\frac{1}{2}(\mathcal{E}_{n,\lambda}+(-1)^n\mathcal{E}_{n,-\lambda})=\left\{
\begin{split}
&1 \quad {\rm{if}} \quad n=0 \\
&0 \quad {\rm{if}} \quad n>0.
\end{split}\right.
\end{equation}

For $n\geq1$, we have
\begin{equation}\label{eq41}
\begin{split}
\mathcal{E}_{n,\lambda}=(-1)^{n-1}\mathcal{E}_{n,-\lambda}.
\end{split}
\end{equation}

From \eqref{eq37} and \eqref{eq41}, we obtain the following theorem. 

\begin{theorem}
Let $a \in \frac{1}{\lambda} E$,\,and let $|\lambda|_{p} \ge 1$. Then we have the following relation.
\begin{equation*}
\begin{split}
-\tanh_\lambda\bigg(\frac{1}{2}:a\bigg)=\sum_{n=1}^\infty\mathcal{E}_{n,\lambda}\frac{a^n}{n!}.
\end{split}
\end{equation*}
\end{theorem}

\medskip

\section{Conclusion}
We intorduced several degenerate hyperbolic functions which are degenerate versions of the usual hyperbolic functions. We computed Volkenborn and the fermionic $p$-adic integrals of the degenerate hyperbolic cosine and the degenerate hyperbolic sine functions. From those results, we were able to derive some identities regarding the degenerate hyperbolic tangent and the degenerate hyperbolic cotangent functions. \par
In recent years, various kinds of tools, like generating functions, combinatorial methods, $p$-adic analysis, umbral calculus, differential equations, probability theory, special functions, analytic number theory and operator theory, have been used in studying special numbers and polynomials, and degenerate versions of them. \par
It is one of our future research projects to continue to explore many special numbers and polynomials and their applications to physics, science and engineering as well as to mathematics.

\hspace{6cm}

%

\noindent{\bf{Acknowledgments}} \\
All authors thank Jangjeon Institute for Mathematical Science for the support of this research.

\vspace{0.1in}

\noindent{\bf {Availability of data and material}} \\
Not applicable.

\vspace{0.1in}

\noindent{\bf{Funding}} \\
The second author is supported by the Basic Science Research Program, the National
               Research Foundation of Korea,
               (NRF-2021R1F1A1050151).
\vspace{0.1in}

\noindent{\bf{Ethics approval and consent to participate}} \\
All authors declare that there is no ethical problem in the production of this paper.

%
%
%
%

%

\


\bigskip
\begin{thebibliography}{9}


\bibitem{1} Araci, S.; Acikgoz, M.; Park, K.-H.; Jolany, H. On the unification of two families of multiple twisted type polynomials by using $p$-adic $q$-integral at $q=-1$. Bull. Malays. Math. Sci. Soc. (2) 37 (2014), no. 2, 543-554.

\bibitem{2} Araci, S.; Erdal, D.; Seo, J. J. A study on the fermionic $p$-adic $q$-integral representation on $\mathbb{Z}_p$ associated with weighted $q$-Bernstein and $q$-Genocchi polynomials. Abstr. Appl. Anal. 2011, Art. ID 649248, 10 pp.

\bibitem{3} Carlitz, L. Degenerate Stirling, Bernoulli and Eulerian numbers. Utilitas Math. 15 (1979), 51-88.

\bibitem{4} Dolgy, D. V.; Park, J.-W. A note on a sum of powers of $q$-integers of skip count by k. Adv. Stud. Contemp. Math. (Kyungshang) 33 (2023), no. 1, 87-93.

\bibitem{5} Kim, D. S.; Kim, T. Some $p$-adic integrals on $\mathbb{Z}_p$ associated with trigonometric functions. Russ. J. Math. Phys. 25 (2018), no. 3, 300-308.

\bibitem{6} Kim, D. S.; Kim, T.; Kwon, J.; Lee, S.-H. Park, S. On $\lambda$-linear functionals arising from $p$-adic integrals on $\mathbb{Z}_p$. Adv. Difference Equ. 2021, Paper No. 479, 12 pp.

\bibitem{7} Kim, T. On the analogs of Euler numbers and polynomials associated with $p$-adic $q$-integral on $\mathbb{Z}_p$ at $q=-1$. J. Math. Anal. Appl. 331 (2007), no. 2, 779-792.

\bibitem{8} Kim, T. $q$-Euler numbers and polynomials associated with $p$-adic $q$-integrals. J. Nonlinear Math. Phys. 14 (2007), no. 1, 15-27.

\bibitem{9} Kim, T.; Kim, D. S. Some identities on degenerate $r$-Stirling numbers via Boson operators. Russ. J. Math. Phys. 29 (2022), no. 4, 508-517.

\bibitem{10} Kim, T.; Kim, D. S. Some identities involving degenerate stirling numbers associated with several degenerate polynomials and numbers. Russ. J. Math. Phys. 30 (2023), no. 1, 62-75.

\bibitem{11} Kim, T.; Kim, D. S.; Park, J.-W. Fully degenerate Bernoulli numbers and polynomials. Demonstr. Math. 55 (2022), no. 1, 604-614.

\bibitem{12} Kwon, J.; Kim, W. J.; Rim, S.-H. On the some identities of the type 2 Daehee and Changhee polynomials arising from $p$-adic integrals on $\mathbb{Z}_p$. Proc. Jangjeon Math. Soc. 22 (2019), no. 3, 487-497.

\bibitem{13} Schikhof, W. H. Ultrametric calculus. An introduction to $p$-adic analysis. Cambridge Studies in Advanced Mathematics, 4. Cambridge University Press, Cambridge, 1984. 

\bibitem{14} Shiratani, K. On some operators for $p$-adic uniformly differentiable functions. Japan. J. Math. (N.S.) 2 (1976), no. 2, 343-353.

\bibitem{15} Shiratani, K.; Yokoyama, S. An application of $p$-adic convolutions. Mem. Fac. Sci. Kyushu Univ. Ser. A 36 (1982), no. 1, 73-83.

\bibitem{16} Woodcock, C. F. An invariant $p$-adic integral on $\mathbb{Z}_p$. J. London Math. Soc. (2) 8 (1974), 731-734.

\bibitem{17} Woodcock, C. F. Fourier analysis for $p$-adic Lipschitz functions. J. London Math. Soc. (2) 7 (1974), 681-693.

\bibitem{18} Yun, S. J.; Park, J.-W. On fully degenerate Daehee numbers and polynomials of the second kind. J. Math. 2020, Art. ID 7893498, 9 pp.







\end{thebibliography}
\end{document}